\documentclass[graybox]{svmult}
\usepackage{comment}

\usepackage{type1cm}        
\usepackage{makeidx}         
\usepackage{graphicx}        
\usepackage{multicol}        
\usepackage[bottom]{footmisc}

\usepackage{newtxtext}       %
\usepackage{newtxmath}       
\usepackage{dsfont}
\usepackage{amsmath}
\usepackage{subcaption}
\usepackage{float}

\makeindex             

\allowdisplaybreaks

\begin{document}

\title*{Neural Drift Estimation for Ergodic Diffusions: Non-parametric Analysis and Numerical Exploration}
\titlerunning{Neural Drift Estimation for Ergodic Diffusions}
\author{Simone Di Gregorio and Francesco Iafrate}
\institute{Simone Di Gregorio \at Department of Computer, Control and Management Engineering Antonio Ruberti,
Sapienza University of Rome, Italy, \email{simone.digregorio@uniroma1.it}
\and Francesco Iafrate \at Department of Mathematics, University of Hamburg, Bundestr. 55 20146 Hamburg, Germany, \email{francesco.iafrate@uni-hamburg.de}}
%
%
\index{Di Gregorio S.} 
\index{Iafrate F.} 

\maketitle

\abstract{We take into consideration generalization bounds in \cite{oga_koike} for the problem of the estimation of the drift component for ergodic stochastic differential equations, when the estimator is a ReLU neural network and the estimation is non-parametric with respect to the statistical model. We show a practical way to enforce the theoretical estimation procedure, enabling inference on noisy and rough functional data. Results are shown for a simulated Itô-Taylor approximation of the sample paths.}

\section{Neural Non-parametric Drift Estimation for Ergodic Diffusions}

The problem of statistical modeling for multi-dimensional ergodic diffusions from discrete observations has a long-standing literature that largely focused on parametric approaches to the estimation of the coefficient set \cite{kessler, iafrate}. In this work, we pick a different approach stemming from the analysis in \cite{oga_koike}, focusing on neural network estimators of the drift coefficient and without parametric assumptions on the statistical model. This nonparametric approach is relevant because, under mild regularity conditions, it enables practical and flexible estimation for an underlying functional variable with a non-trivial noise dynamic, from which we sample points. Moreover, the generality, rough path features and analytical tractability of diffusions make them a standard tool for modeling structured and irregular time sequences as random functions. While finance is the standard application, as stochastic calculus provides its mathematical framework, these methods are also valuable for investigating data in broader contexts, such as biological or geological data analysis \cite{functional_data_bio, source_function_sde}.

The formal setting is the following. Let $\{W_t\}_{t\in[0, T]}$ be a $d$-dimensional Brownian motion adapted to an underlying filtered probability space, 
$\sigma: \mathbb{R}^d\rightarrow \mathbb{R}^{d\times m}$ and $b:\mathbb{R}^d \rightarrow \mathbb{R}^d$, which for our purposes we assume as globally Lipschitz (\textbf{B1}):

\begin{equation}
\label{eq::sde_definition}
   dX_t = b(X_t)dt + \sigma(X_t)dW_t, \, X_0 = \eta
\end{equation}
According to standard results, under \textbf{B1}, a strong solution $\{X_t\}_{t\in[0, T]}$ to the above stochastic differential equation (SDE) always exists. We also require an exponential mixing condition (\textbf{B2}): 
$\exists\,  C_\beta, C'_\beta: \beta_X(t)\leq C'_\beta\exp\left(-C_\beta t\right)$. $\beta_X$ is the $\beta$-strong mixing function of the solution. A useful sufficient condition \cite{beta_mixing_decay} for \textbf{B2} is that we require that $\eta$ has a moment generating function, that $\sigma\sigma^{\mathsf{T}}$ is bounded and bounded away from zero and that $x^\mathsf{T}b(x)\leq -r \lVert x\rVert_2^{\alpha} \, \forall \, x \, \text{s.t.}\, \lVert x\rVert_2 \geq M_0$, with $M_0 \geq 0, \ \alpha \geq 1, \ r>0$.

The estimation target is $b$, with $\sigma$ not entering the analysis. The estimation is restricted to a parallelotope, which can be taken without loss of generality as $[0, 1]^d$: one may pick a generic parallelotope and scale and shift accordingly. The target is thus a specific restricted component of the drift, $f_0 := b_i \mathds{1}_{[0, 1]^d}$, for $i \in \{1, \dots, d\}$.

Following \cite{oga_koike}, we formalize a neural network as a tuple $(L, \mathbf{p}, \{W_i\}_{i=0}^L, \{\mathbf{v}_j\}_{j=1}^L)$, and the non-linearity is taken to be the \textit{Rectified Linear Unit} (ReLU), i.e. $\rho_{v_i}(x_i):= \max(0, x-v_i)$, where $v_i$ is the $i$-th entry of a generic shift vector $\mathbf{v}$. Every tuple like the above induces a map $f$ in the following standard way:
\begin{align}
\label{eq::neural_network_definition}
    f(x) &:= W_L\rho_{\mathbf{v}_L}W_{L-1}\rho_{\mathbf{v}_{L-1}} \dots W_{1}\rho_{\mathbf{v}_{1}}W_{0}x, \nonumber \\ 
    & x\in \mathbb{R}^d, \,  W_l \in \mathbb{R}^{p_{l+1}\times p_l}, \, p_0 = d,\, p_{L+1} = 1 \nonumber \\
    \rho_{\mathbf{v}}(y) &:= \begin{pmatrix}
        \rho_{v_1}(y_1) \\
        \vdots \\
        \rho_{v_r}(y_r)        
    \end{pmatrix}, \ y\in \mathbb{R}^r
\end{align}

The hypothesis class of the estimator we consider is: 
\begin{align}
\label{eq::hypothesis_class}
    \mathfrak{F}(L, \mathbf{p}, s, F) :=& \Bigg\{f\mathds{1}_{[0, 1]^d}: f\text{ follows  (\ref{eq::neural_network_definition}), }\underset{j=0, \dots, L}{\max} \lVert W_j\rVert_\infty \lor \lVert \mathbf{v}_j\rVert_\infty\leq 1, \nonumber \\
    & \lVert W_0\rVert_0 + \sum_{j=1}^L (\lVert W_j\rVert_0 + \lVert \mathbf{v}_j\rVert_0)\leq s, \lVert f\mathds{1}_{[0, 1]^d}\rVert_\infty\leq F\Bigg\}
\end{align}

Given a constant sampling interval $\delta$ and a number of points $N$, we observe a sample of the solution $\{X_{k\delta}\}_{k=0}^N$. The idea is defining a response variable for learning $f_0$ as $Y_{k\delta}:= \frac{1}{\delta}\left(X^i_{(k+1)\delta}-X^i_{k\delta}\right)$. The almost sure rough path properties of the SDE solution and the exploding variance of the difference quotients $Y_{k\delta}$ present a challenge in this context, see the analysis in Subsection \ref{subsec::analysis_results}. One may connect this reasoning to the standard statistical learning scenario, showing that the drift function in $X_{k\delta}$ can be seen as an increasingly precise (as $\delta \downarrow 0$) approximation of the regression function for the prediction of $Y_{k\delta}$. This in turn heuristically justifies a mean square empirical risk minimization, which targets $f_0$ through the proxy $Y_{k\delta}$. This leads to a training loss and a generalization risk for the estimator $\hat{f}_N$ given by:
\begin{align}
\label{eq::loss_definition}
    \mathcal{Q}_N(f) := & \frac{1}{N}\sum_{k=0}^{N-1}\left(Y_{k\delta}-f(X_{k\delta})\right)^2 \\
\label{eq::risk_definition}
    \mathcal{R}_N(\hat{f}_N, f_0):= & \,\mathbb{E}\left[\frac{1}{N}\sum_{k=0}^{N-1}\left(\hat{f}_N(X'_{k\delta})-f_0(X'_{k\delta})\right)^2\right]
\end{align}
where $\{X'_t\}_{t\in [0, N\delta]}$ is an independent copy  of the solution to  (\ref{eq::sde_definition}).

In \cite{oga_koike} the authors prove that under \textbf{B1}, \textbf{B2}, with $F \geq \lVert f_0\rVert_\infty$, any estimator taking values in $\mathfrak{F}(L, \mathbf{p}, s, F)$ satisfies the following risk inequality, for all $\tau > 0$:
\begin{align}
\label{eq::risk_guarantee}
\mathcal{R}_N(\hat{f}_N, f_0)\leq &\,\tau\left(\Psi_N^{\mathfrak{F}(L, \mathbf{p}, s, F)}(\hat{f}_N) + \underset{f \in \mathfrak{F}(L, \mathbf{p}, s, F)}{\inf}\lVert f - f_0\rVert\right) + \nonumber \\
&+ C_\tau F^2\left(\frac{s(L\log{s} + \log{(N\delta)})\cdot \log{(N\delta)}}{N\delta} + \delta\right)   
\end{align}
where in the above $\Psi_N^{\mathfrak{F}(L, \mathbf{p}, s, F)}$ is the expected optimization error of the fitting procedure with respect to the considered class $\mathfrak{F}$, while $C_\tau$ is a constant depending on $\tau$, the specification of \textbf{B1} and \textbf{B2} and the $\lVert\cdot\rVert_\infty$ norm of the coefficients over the estimation set. The last term - which represents the estimation variance - vanishes under the standard asymptotic regime \cite{kessler}, i.e. $ T = N\delta \uparrow \infty$ and $\delta \downarrow 0$.

\section{Building Regular Networks}
\label{sec::2}

The problem of enforcing the regularity specified in  (\ref{eq::hypothesis_class}) is relevant. Clearly, we could just impose an appropriate hard norm constraint on blocks of the weight set and take $s$ as the total cardinality of the parameter set. For example, for a specific parameter vector, constraining the estimate inside a $\ell^2$ ball of unitary radius would suffice to enforce each entry being bounded by 1. This clearly strongly limits the flexibility of the model leading to increased depth or width to reduce the approximation error. We aim for a much more compact representation of the networks while preserving the same regularity as the one needed for the risk inequality in  (\ref{eq::risk_guarantee}) to hold. The following lemma is Proposition A.3 from \cite{deep_neural_approximation}, with an additional remark from \cite{oga_koike}:

\begin{lemma}
\label{lemma::network_conversion}
    For any ReLU neural network function $g:\mathbb{R}^d \rightarrow \mathbb{R}$ with parameters bounded by $B$, there exists another ReLU network $\phi$ such that $\phi(x) = g(x)\, \forall\, x \in \mathbb{R}^d$ and such that the weights of $\phi$ are uniformly bounded by 1. Additionally, with $L$, $s$ and $p$ respectively the depth, the number of non-null parameters and the maximum width of $g$, the depth of the new network $\phi$ is $L'=\lceil (\log B + 5)L\rceil$, its width $p'$ is bounded by $\max\{3, p\}$ and the number of active parameters is at most $s'= 2s + 12L'$.
\end{lemma}

Therefore, we get a logarithmic growth in $B$ for all the complexity metrics that appear in  (\ref{eq::risk_guarantee}). What this suggests is that if we are able to enforce any uniform bound on the weight set we can still get that the produced estimator $\hat{f}_N\in \mathfrak{F}(L', \mathbf{p}', s', F)$, just by considering the regular version of the estimator stemming from Lemma \ref{lemma::network_conversion}.
The next Lemma shows that a bounding constant $F$ always exists as soon as we impose, while optimizing, a norm constraint on each element of the parameter set.

\begin{lemma}
\label{lemma::equiboundness}
    Let $K$ be any set of indexes, let $\{(L, \mathbf{p}, \{W_{ik}\}_{i=0}^L, \{\mathbf{v}_{jk}\}_{j=1}^L)\}_{k\in K}$ be a set of neural network tuples with input in $\mathbb{R}^d$ for which $\exists D>0: \lVert W_{ik} \rVert \leq D, \lVert \mathbf{v}_{jk} \rVert \leq D \,\forall \,i \in\{0, \dots, L\}, j \in\{1, \dots, L\}, k\in K$. Let $\{g_k\}_{k\in K}$ be the set of functions which are induced by each of the networks above. Then we have that $\exists F\geq 0: \lVert g_k \mathds{1}_A\rVert_\infty \leq F \, \forall \, k \in K$, where $A \subset \mathbb{R}^d$ and is compact.
\end{lemma}
The result can be proved with induction over composition and sub-multiplicativity of matrix norms. By combining Lemma \ref{lemma::network_conversion} and Lemma \ref{lemma::equiboundness} we get the following Theorem: 
\begin{theorem}
    Let the underlying diffusion process satisfy \textbf{B1} and \textbf{B2}, and let $\hat{f}_N$ be a feedforward ReLU neural network estimator with weights bounded by an arbitrary constant $B$. Then the estimator is bounded on $[0, 1]^d$ by a constant $F$ whose existence is assured through Lemma \ref{lemma::equiboundness} and thus satisfies - if $F \geq\lVert f_0 \rVert_\infty$ - the risk inequality in  (\ref{eq::risk_guarantee}), with the configuration specified by Lemma \ref{lemma::network_conversion}.
\end{theorem}

Notice that the optimization error and the approximation error in  (\ref{eq::risk_guarantee}) need to be taken w.r.t. $\mathfrak{F}(L', \mathbf{p}', s', F)$. These results enable much more flexibility when tuning the optimization procedure: for example, in our case, setting $L=2$ and $32$ as hidden dimensionality, we enforce the required regularity conditions just by imposing $\lVert \mathbf{v}_l\rVert_2\leq \sqrt{p_l}$  for $l \in \{1, 2\}$ and $\lVert W_{l_{j, \cdot}}\rVert_2\leq \sqrt{p_l}$ for $\, j\in \{1, \dots, p_{l+1}\},\, l \in \{0, 1, 2\}$, which is the $\ell^2$ bound we would immediately get if all the entries of each row in each weight component were bounded by 1. This still clearly implies that the $\lVert \cdot \rVert_\infty$ norms in  (\ref{eq::hypothesis_class}) are bounded, but by $\max_{l\in\{0, \dots, L\}}\sqrt{p_l}$ and not by 1. The approach of constraining the $\ell^2$ norm of the weight vector of each linear combination is a common heuristic idea when training neural networks \cite{dropout}, being used to drive regularization on the estimator. 
The $\ell^2$ constraint is easily enforced when training with stochastic gradient descent on the loss from  (\ref{eq::loss_definition}), simply by rescaling after each update if needed.

\textit{Further developments}. Other approaches to neural network regularization with statistical guarantees have been studied e.g. in \cite{statistical-guarantees}. Future developments of our work can be an extension of such results to our setting. 

\section{Numerical Monte Carlo Testing}

In order to check the validity of the theoretical results we introduce a testing framework based on the SDE model:

\begingroup\makeatletter\def\f@size{9}
\begin{equation}
\label{eq::numerical_test_sde_rescaled}
dY_t = \left[ 
\begin{matrix}
-\alpha_1 Y_{t, 1} + c_1\alpha_2 \left(\sin\left(\frac{Y_{t, 2}}{c_2}\right) + 2\right) \\
c_2\alpha_3 \left(\cos\left(\frac{Y_{t, 1}}{c_1}\right) + 2\right) - \alpha_4 Y_{t, 2}
\end{matrix} \right] dt\,
+\left[
\begin{matrix}
\beta_1 c_1s\left(\frac{Y_{t, 1}}{c_1}\right) + c_1\beta_3 & 0 \\
0 & \beta_2 c_2s\left(\frac{Y_{t, 2}}{c_2}\right) + c_2\beta_3
\end{matrix} 
\right] dW_t
\end{equation}
\endgroup

This model is an extension of the one proposed in \cite{yoshida-2012} to the case of non-constant diffusion.
In the above, $\alpha_1 = 1, \alpha_2 = 2, \alpha_3 = 2, \alpha_4 = 1, \beta_1 = 0.5, \beta_2 = 0.5, \beta_3 = 0.1$ are hyperparameters that specify the probabilistic behavior of the solution; specifically, it is crucial that $\alpha_1$ and $\alpha_4$ are negative and that the diagonal of the diffusion coefficient remains bounded away from zero, something that is achieved with the $\beta_3$ offset. $c_1 = \frac{1}{6}$ and $c_2 = \frac{1}{5}$ are instead much less relevant and, as anticipated before, are just constants that results from scaling the SDE solution for it to reside with high probability mostly in $[0, 1]^2$.
\label{sec::numerical_testing}
\subsection{Numerical Setup}
\label{subsec::numerical_setup}

With a deterministic initial condition (in our case $Y_0 = 0$) and our hyperparameter choice, it can be shown that  (\ref{eq::numerical_test_sde_rescaled}) satisfies the needed ergodic conditions from \cite{beta_mixing_decay}; additionally, this SDE model satisfies the linear growth and Lipschitz conditions on the first Itô-Taylor coefficients \cite{kloeden_ito}, which in turn allow an approximation of the sample paths that is $L^1$ convergent to the solution, the so-called Milstein scheme.

We fix a time approximation mesh $\Delta=10^{-3}$. The sampling of the solution is induced by skipping points evenly in the time discretization. This allows us to mitigate the spurious error induced by the numerical computation of the sample paths of the solution to  (\ref{eq::numerical_test_sde_rescaled}). Having fixed $\Delta$, the sampling $\delta$ appearing in (\ref{eq::risk_guarantee}) is determined by the number of skipped points, which we call \textit{skip} in the following. 

The actual implementation of the neural network, according to the configuration specified in Section \ref{sec::2}, is coded in Python with the Tensorflow framework. The number of epochs for the (mini-batch) Adam stochastic gradient descent routine is set fixed to 200, while the number of Monte Carlo iterations is set fixed to 50. After each training routine and evaluation on the test trajectory, we save the error metric contained in the expectation in  (\ref{eq::risk_definition}) for Monte Carlo approximation.

Notice also that we build a single network to estimate both components together, something that differs from the treatment in \cite{oga_koike}: the idea is that we can split the network in two by taking slices of the last dense layer. This induces two networks that both satisfy the regularity we have been describing, sharing the same parameter set except for the last linear map. 

\subsection{Analysis of Results}
\label{subsec::analysis_results}
\begin{figure}[H]
    \centering
    \begin{minipage}{0.45\textwidth}
        \centering
        \begin{subfigure}[b]{\textwidth}
            \centering
            \includegraphics[width=\textwidth]{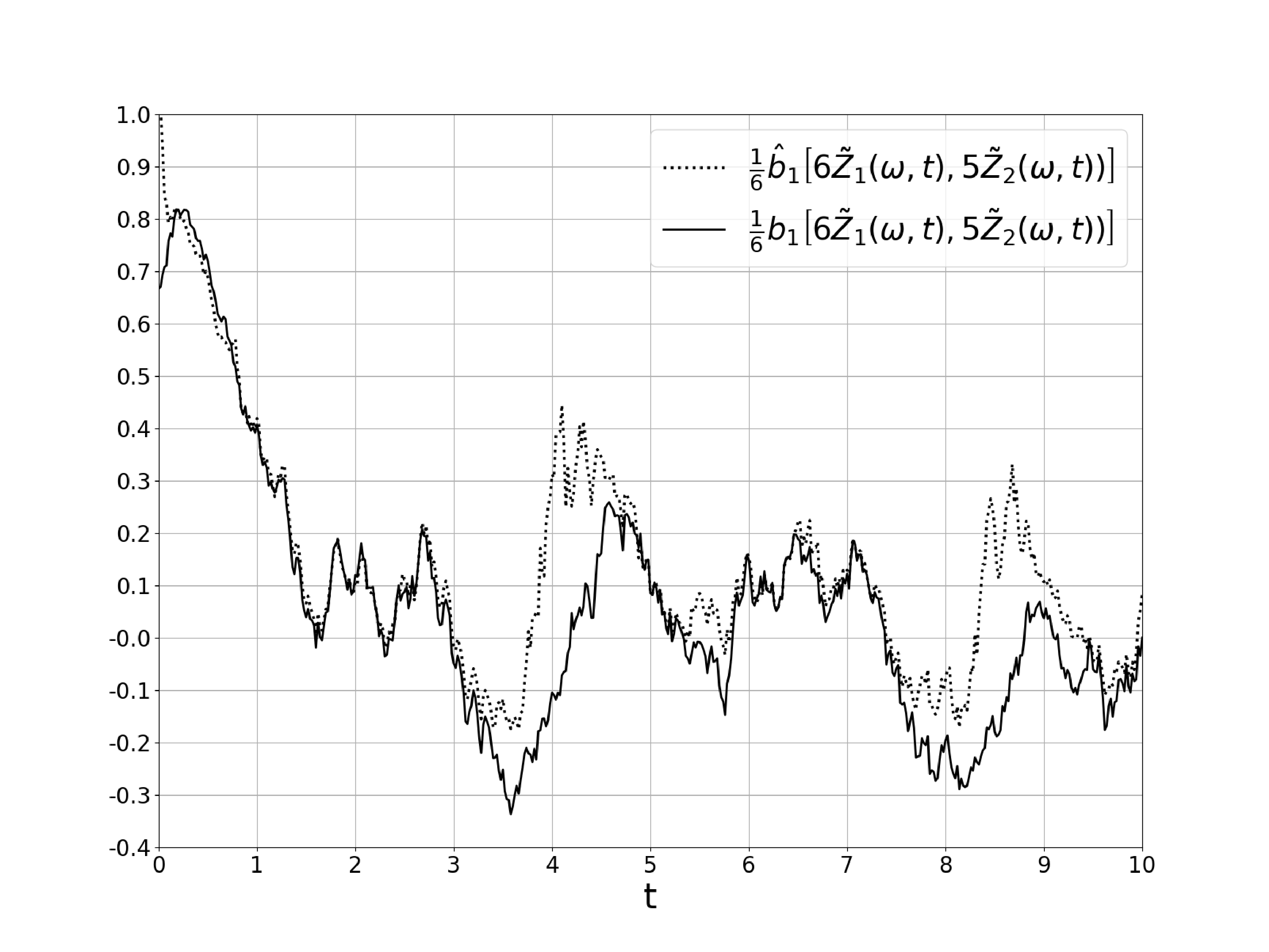}
            \caption*{First component of the drift, $T=10$}
            \label{fig::test_drift_1_10}
        \end{subfigure}
        
        \begin{subfigure}[b]{\textwidth}
            \centering
            \includegraphics[width=\textwidth]{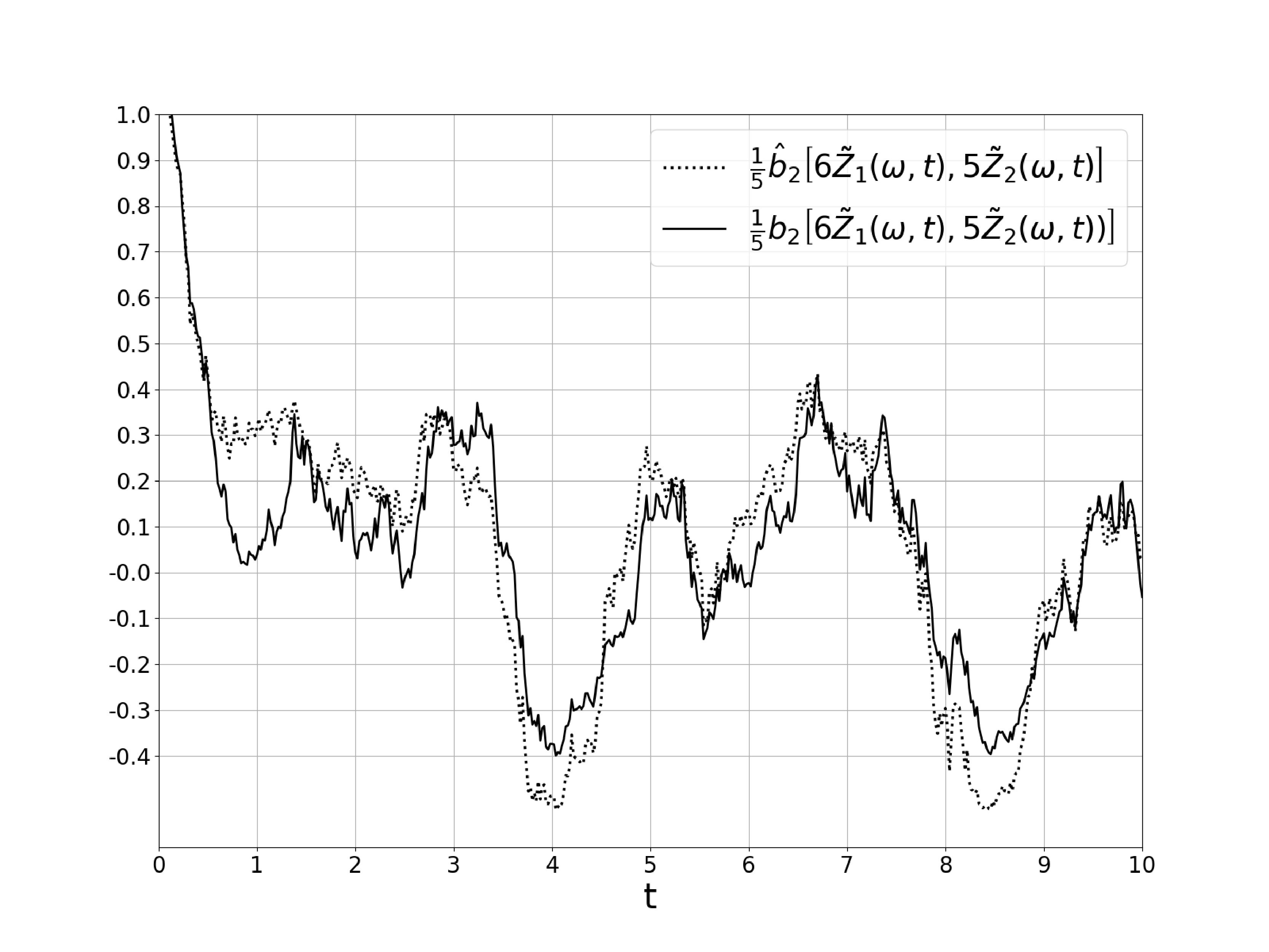}
            \caption*{Second component of the drift, $T=10$}
        \end{subfigure}
        
        \caption*{\textbf{(a)} Drift trajectories for $T=10$}
    \end{minipage}
    \hfill
    \begin{minipage}{0.45\textwidth}
        \centering
        \begin{subfigure}[b]{\textwidth}
            \centering
            \includegraphics[width=\textwidth]{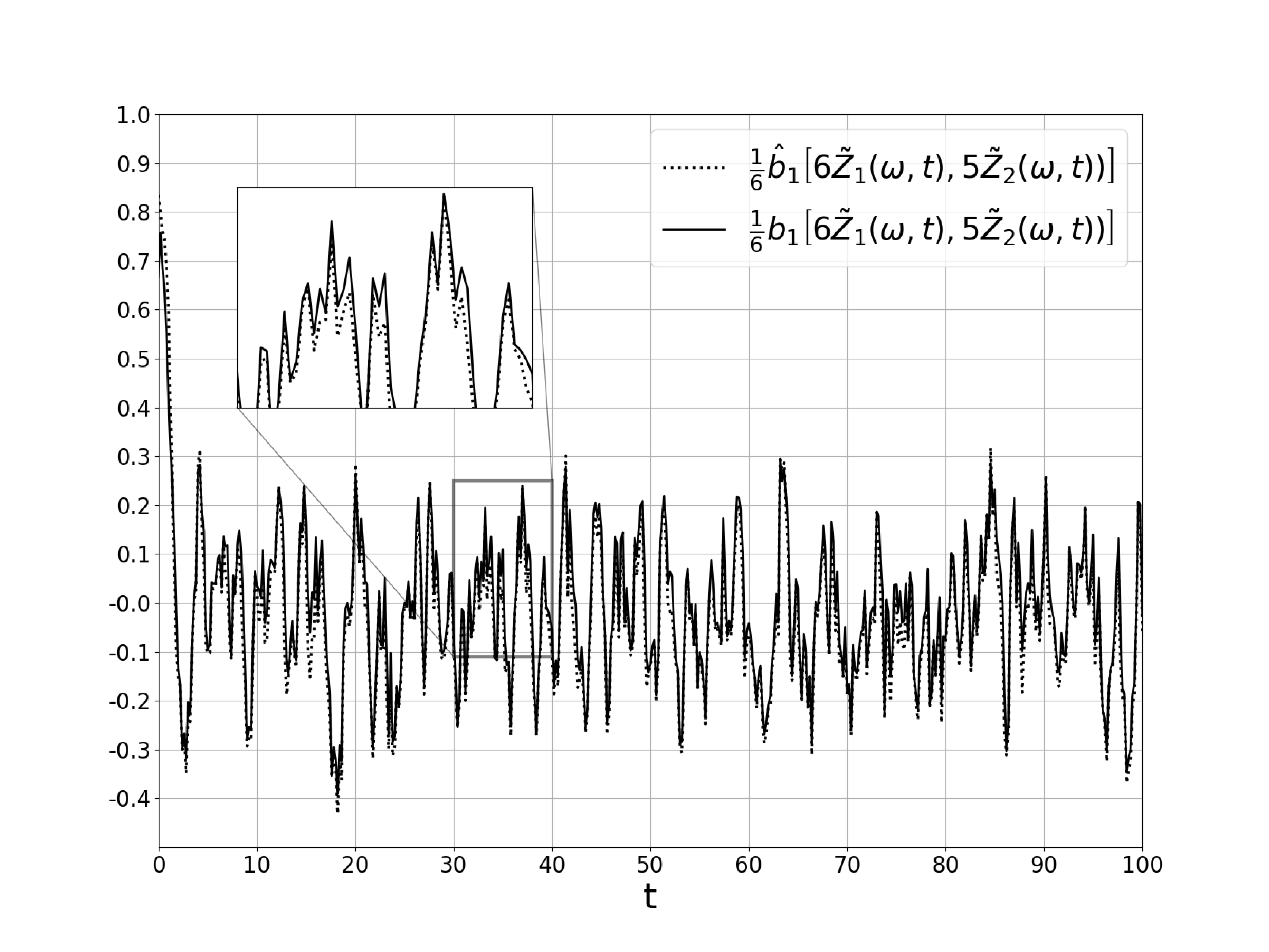}
            \caption*{First component of the drift, $T=100$}
        \end{subfigure}
        
        \begin{subfigure}[b]{\textwidth}
            \centering
            \includegraphics[width=\textwidth]{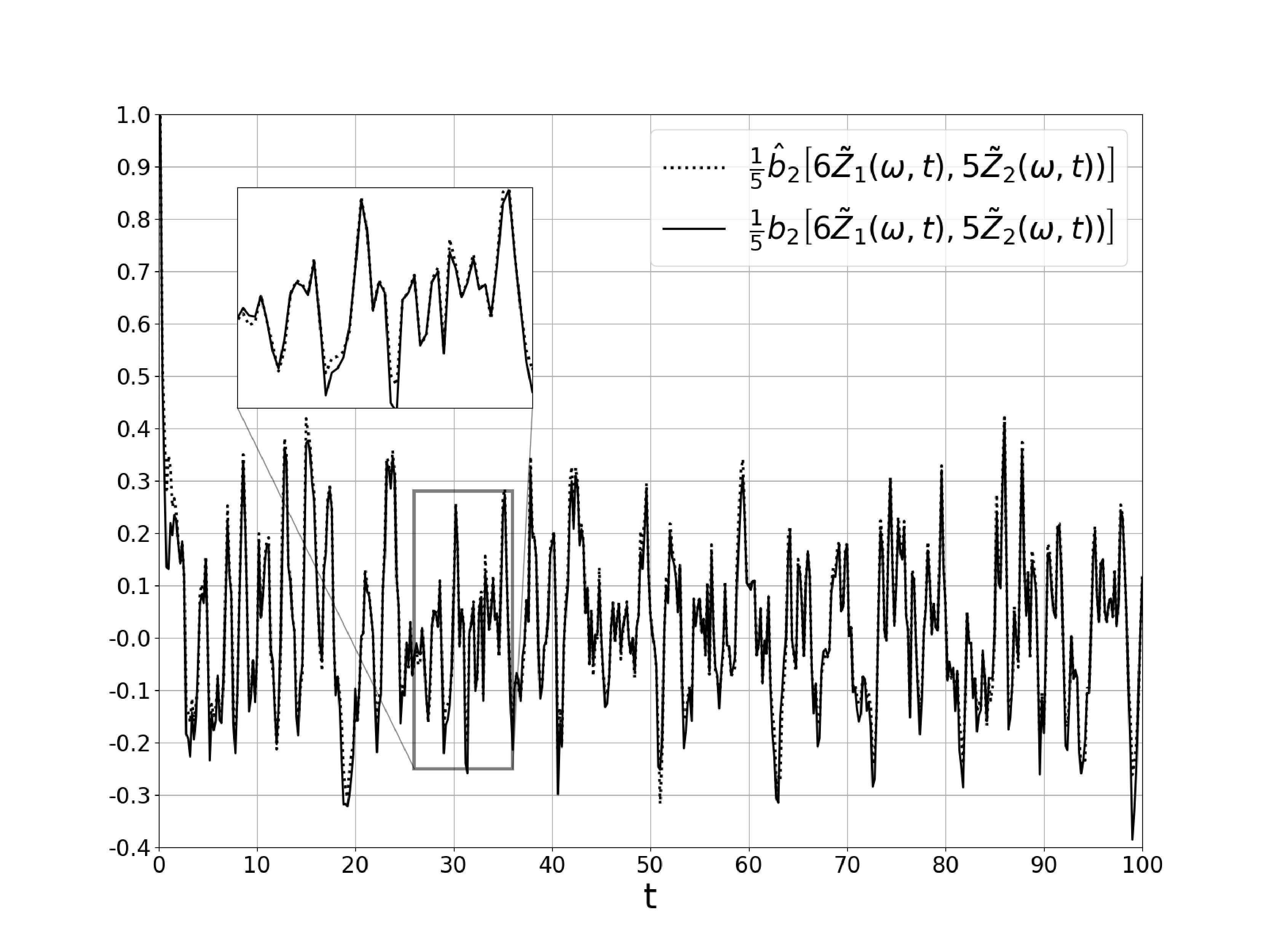}
            \caption*{Second component of the drift, $T=100$}
        \end{subfigure}
        
        \caption*{\textbf{(b)} Drift trajectories for $T=100$}
    \end{minipage}
    \captionsetup{labelfont=bf}
    \caption{Actual and modeled trajectories for the drift components evaluated on a sample path of the approximation of the test process, $\{\tilde{Z}_t\}_t$; $20$ skipped observations and $T=10$ (left) and $T=100$ (right). $b_1$ and $b_2$ in the plot legends are the components from the non-scaled SDE (i.e. $c_1 = c_2 = 1$) from  (\ref{eq::risk_guarantee})}
    \label{fig::sample_paths_test}
\end{figure}

In Fig. \ref{fig::sample_paths_test}, 
we show a single Monte Carlo iteration and get some visual insights.
We train and test with two different time horizons $T=10, 100$, $\delta = 0.02$, corresponding to \textit{skip} $= 20$.   
The impact of $T$ on improving the fit arises. 
This can also be seen from the final results in Table \ref{tab::mse}.
This is expected, since the increase in $T$ leads to the exploration of the high-probability sets under the ergodic measure of the process. For low $T$, ergodic behavior cannot help much, although Fig. \ref{fig::sample_paths_test}.a shows that a crude modeling of the relevant section of the drift surface can still be achieved. This is also not surprising due to the exponentially fast total variation convergence of the Markov kernels to the ergodic measure, which is again assured here \cite{beta_mixing_decay}.

\begin{figure}[H]
    \centering
    \captionsetup{labelfont=bf}
    \begin{subfigure}[t]{0.48\textwidth} 
        \centering
        \includegraphics[width=\textwidth]{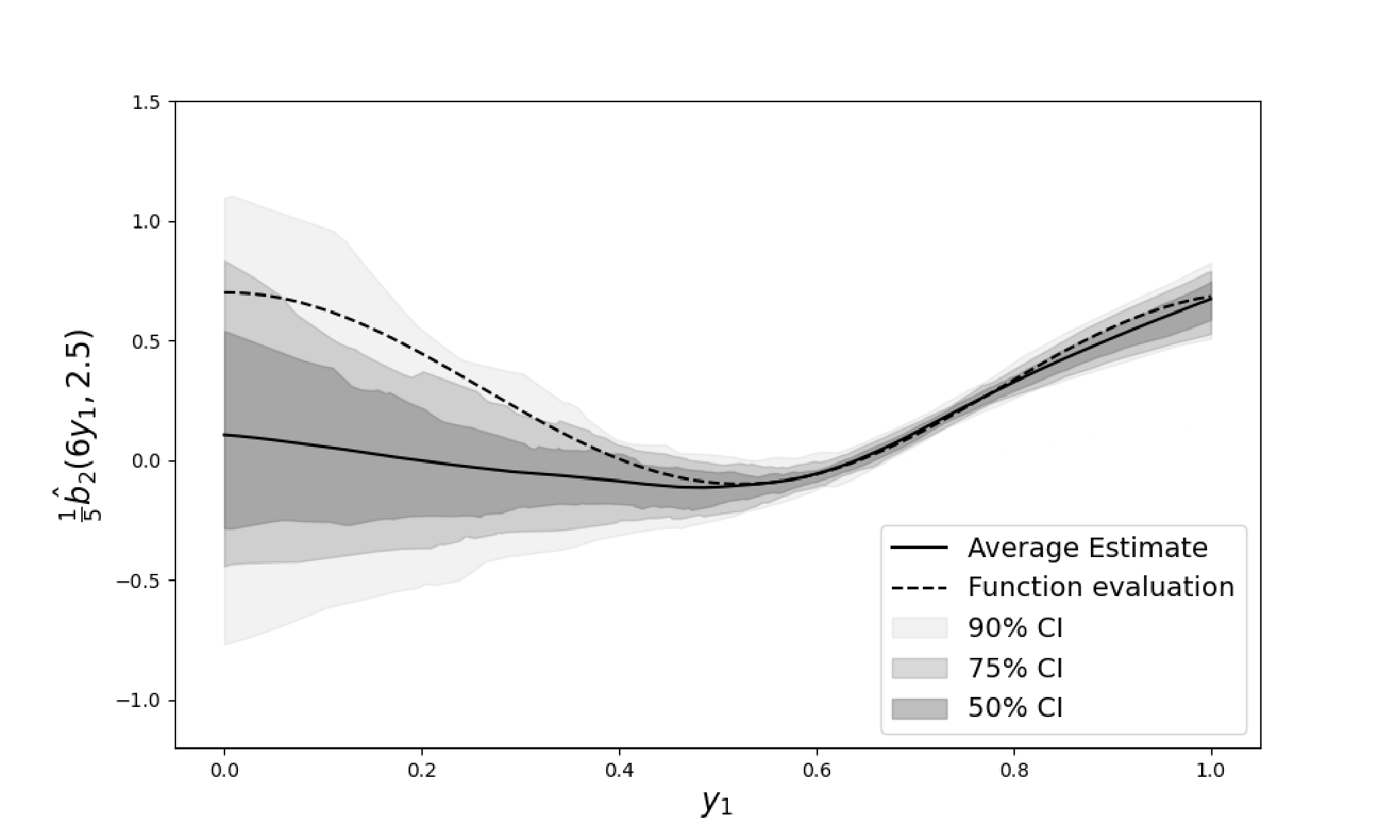}
        \caption{Average estimate, confidence bands and actual evaluation for the second component of the drift, when fixing the second of the two input variables}
        \label{fig::drift_cut_predicted}
    \end{subfigure}
    \hfill
    \begin{subfigure}[t]{0.48\textwidth} 
        \centering
        \includegraphics[width=\textwidth]{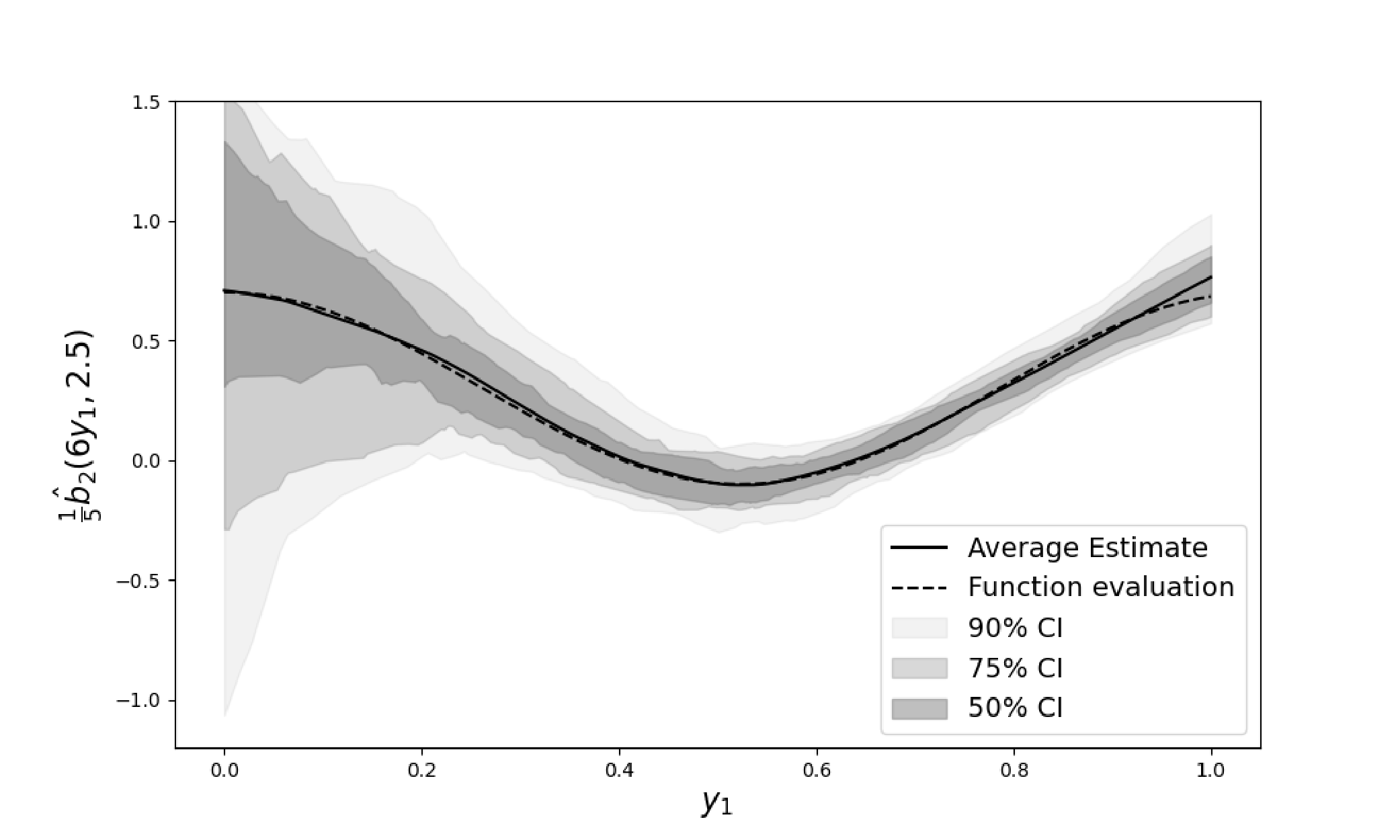}
        \caption{Average estimate, confidence bands and actual evaluation for the second component of the drift, when fixing the second input variable and with diffusion scaled by $2.5$}
        \label{fig::drift_cut_predicted_noise}
    \end{subfigure}
    \caption{A 50 iterations Monte Carlo experiment: we compare the function surface with the estimated one, across two SDE configurations where for one we are scaling the diffusion component. For the fit, the configuration is the one of Fig. \ref{fig::sample_paths_test}.b}
    \label{fig::estimation_bias}
    
\end{figure}

In the interpretation of Fig. \ref{fig::sample_paths_test} it is important to note that the target of our estimate is the drift evaluated on the random trajectory, not the drift itself uniformly over the $[0, 1]^2$ domain. This means that we cannot expect to be able to estimate state space sets that are visited with low probability. This immediate intuition is actually shown in practice with another experiment, which is highlighted in Fig. \ref{fig::estimation_bias}. The idea is that upward scaling the diffusion coefficient should increase the noise component of the process, leading it to a more uniform coverage of the considered box domain. 

In Fig. \ref{fig::estimation_bias} we show
the variance and the bias of (a slice of) the estimated surface $\hat b_j(y_1, y_2), j=1,2$. The Monte Carlo simulation supports this idea of the drift being learned non-homogeneously, with a strong estimation bias in some regions. This effect is indeed mitigated when increasing the magnitude of the noise stemming from the diffusion component.

\vspace{-1.2cm}

Table \ref{tab::mse} gives an overall view of the performance in terms of the risks involved. Apart from the obvious relevance of $T$, the analysis highlights a counter-intuitive behavior with the decreasing $\delta$ (i.e. decreasing \textit{skip} variable). Specifically, \textit{Test MSE} seems not to benefit from low values of $\delta$, after decreasing the \textit{skip} variable past $100$.

This would seem in contrast with the finite-sample valid bound in  (\ref{eq::risk_guarantee}). Nevertheless, we stress that the bound \textbf{may not be tight} in our setting due to unspecified universal constants. This non-tightness may be caused by the model being moderately parameterized. This in turn may cause overfitting problems arising 
when fitting the difference quotients, since a vanishing $\delta$ increases the variance of the response (recall that the underlying process is not differentiable).
The error on the quotients (last column in Table \ref{tab::mse}) is compatible with the irreducible error from Table \ref{tab::irreducible_noise},
with the discrepancy between the two indeed increasing for \textit{skip} $< 100$, with the fit starting to capture noise.
This explains the non-monotonic in $\delta$ behaviour of \textit{Train MSE}.

It is apparent that the difference between \textit{Train MSE} and \textit{Test MSE} also increases, meaning that the estimated surface is greedily representing the drift evaluated in the points of the sample. Our hypothesis is that this happens since increasing $\delta$ leads to datasets with observations that are close to each other in the state space, with the result similar to increasing the training iterations, by continuity of the drift trajectory.

\begin{table}[H]
\centering
\caption{Average MSE (scaled by $10^3$ and averaged across the two components together, as in Table \ref{tab::mse}) between the difference quotients and the actual drift evaluation, computed with $1000$ Monte Carlo iterations for trajectories with $T=100$}
\label{tab::irreducible_noise}
\begin{tabular}{|c c c c c|}
\hline
\multicolumn{5}{|c|}{Estimate of the irreducible error ($\times 10^{3}$)} \\
\hline
\ skip = 200\  & \ skip = 100 \  & \ skip = 50 \ & \ skip = 20 \ & \ skip = 10 \  \\
\hline
48.297 & 104.348 & 218.657 & 562.421 & 1136.375 \\
\hline
\end{tabular}
\vspace{0.5cm}
\end{table}
\begin{table}[H]
\centering
\caption{The MSEs are computed with 50 Monte Carlo iterations. \textit{Test MSE} is the estimate of the risk  (\ref{eq::risk_definition}), \textit{Train MSE} is the same but with evaluation on the training process and \textit{Quotients MSE} is the error on the training difference quotients from (\ref{eq::loss_definition})}
\label{tab::mse}
\begin{tabular}{|r r r r r|}
\hline
\textbf{skip} & \textbf{T} & \textbf{Test MSE ($\times 10^3$)} & \textbf{Train MSE ($\times 10^3$)} & \textbf{Quotients MSE ($\times 10^3$)} \\
\hline
200 & 10  & \textbf{16.353} (1.204)   & 14.417 (0.629)   & 49.800 (1.281)    \\
200 & 25  & \textbf{8.627} (0.439)    & 7.818 (0.266)    & 48.627 (0.700)    \\
200 & 50  & \textbf{4.557} (0.282)    & 4.177 (0.237)    & 46.461 (0.528)    \\
200 & 100 & \textbf{2.269} (0.107)    & 2.193 (0.093)    & 47.032 (0.279)    \\
\hline
100 & 10  & \textbf{10.995} (0.575)   & 10.148 (0.441)   & 101.165 (1.463)   \\
100 & 25  & \textbf{4.648} (0.353)    & 4.305 (0.195)    & 101.971 (0.870)   \\
100 & 50  & \textbf{2.625} (0.110)    & 2.538 (0.104)    & 101.546 (0.605)   \\
100 & 100 & \textbf{2.013} (0.140)    & 1.786 (0.107)    & 104.299 (0.444)   \\
\hline
50  & 10  & \textbf{11.773} (0.737)   & 10.295 (0.447)   & 215.721 (2.014)   \\
50  & 25  & \textbf{5.190} (0.498)    & 3.900 (0.182)    & 215.832 (1.340)   \\
50  & 50  & \textbf{3.430} (0.324)    & 2.849 (0.180)    & 216.440 (0.851)   \\
50  & 100 & \textbf{2.328} (0.162)    & 2.083 (0.135)    & 218.162 (0.745)   \\
\hline
20  & 10  & \textbf{12.450} (1.155)   & 8.667 (0.416)    & 552.487 (3.776)   \\
20  & 25  & \textbf{6.529} (0.494)    & 5.128 (0.327)    & 552.572 (2.368)   \\
20  & 50  & \textbf{3.706} (0.222)    & 2.838 (0.139)    & 559.533 (1.628)   \\
20  & 100 & \textbf{2.570} (0.232)    & 2.072 (0.155)    & 561.386 (0.914)   \\
\hline
10  & 10  & \textbf{13.894} (1.597)   & 8.527 (0.384)    & 1101.550 (5.237)  \\
10  & 25  & \textbf{7.061} (0.807)    & 4.996 (0.299)    & 1121.530 (2.998)  \\
10  & 50  & \textbf{4.309} (0.430)    & 3.268 (0.256)    & 1129.400 (1.940)  \\
10  & 100 & \textbf{3.706} (0.426)    & 2.375 (0.157)    & 1131.730 (1.327)  \\
\hline
\end{tabular}
\end{table}

\section*{Acknowledgments}
Simone Di Gregorio's research was partly supported by the PNRR MUR project IR0000013-SoBigData.it. This version of the contribution has been accepted for publication, after peer review but is not the Version of Record and does
not reflect post-acceptance improvements, or any corrections. The Version of Record is available online at:
https://doi.org/10.1007/978-3-031-92383-8{\_}20. Use of this Accepted Version is subject to the publisher’s Accepted Manuscript terms of use https://www.springernature.com/gp/open-research/policies/accepted-manuscript-terms.

\end{document}